\documentclass[12pt]{article}
\usepackage[utf8]{inputenc}
\usepackage{amssymb,amsmath,amsfonts,eucal,mathrsfs,amsthm} 
\usepackage[colorinlistoftodos]{todonotes}
\usepackage[allcolors=blue]{hyperref}
\usepackage{empheq}
\newtheorem{theorem}{Theorem}
\newtheorem{proposition}[theorem]{Proposition}
\newtheorem{lemma}[theorem]{Lemma}

\newtheorem{remarks}[theorem]{Remarks}

\theoremstyle{definition}
\newcommand{\R}{\mathbb{R}}
\newcommand{\Q}{\mathbb{Q}}
\newcommand{\Sf}{\mathbb{S}}

\newcommand{\spa}{\mbox{span}}
\newcommand{\hess}{\mbox{Hess\,}}

\newcommand{\grad}{\mbox{grad\,}}

\newcommand{\End}{\mbox{End}}
\newcommand{\Hom}{\mbox{Hom}}
\newcommand{\trace}{\mbox{tr\,}}

\newcommand{\Les}{\mathbb{L}}
\def\<{{\langle}}
\def\>{{\rangle}}

\def\a{\alpha}

\def\id{\mbox{id}}

\def\be{\begin{equation} }
\def\ee{\end{equation} }

\def\proof{\noindent{\it Proof:  }}
\def\qed{\ifhmode\unskip\nobreak\fi\ifmmode\ifinner
\else\hskip5 pt \fi\fi\hbox{\hskip5 pt \vrule width4 pt
height6 pt  depth1.5 pt \hskip 1pt }}
\makeatletter
\newcommand{\subjclass}[2][]{\let\@oldtitle\@title
\gdef\@title{\@oldtitle\footnotetext{#1 
\emph{Mathematics Subject Classification:} #2}}}
\newcommand{\keywords}[1]{\let\@@oldtitle\@title
\gdef\@title{\@@oldtitle\footnotetext
{\emph{Key words and phrases.} #1.}}}
\makeatother
\begin{document}

\title{On the Moebius deformable hypersurfaces}
\maketitle
\begin{center}
\author{M. I. Jimenez        and
        R. Tojeiro$^*$}  
        \footnote{Corresponding author}
      \footnote{This research was initiated while the first author was supported by CAPES-PNPD Grant 88887.469213/2019-00 and was finished under the support of Fapesp Grant 2022/05321-9. The second author was partially supported by Fapesp Grant 2016/23746-6 and CNPq Grant 307016/2021-8.\\
      Data availability statement: Not applicable.}
\end{center}
\date{}

\begin{abstract}
   In the article [\emph{Deformations of hypersurfaces preserving the M\"obius metric and
a reduction theorem}, Adv.  Math. 256 (2014), 156--205], Li, Ma and Wang investigated the interesting class of Moebius deformable hypersurfaces, that is, the umbilic-free Euclidean hypersurfaces $f\colon M^n\to \mathbb{R}^{n+1}$  that admit non-trivial deformations preserving the Moebius metric. The classification of Moebius deformable hypersurfaces of dimension $n\geq 4$ stated in the aforementioned article, however, misses a large class of examples. In this article we complete that classification for $n\geq 5$. 
\end{abstract}

\vspace{2ex}


\date{}
\maketitle

\section{Introduction}
Let $f\colon M^n\to\R^{m}$ be an isometric immersion of a Riemannian manifold $(M^n,g)$ into Euclidean space with normal bundle-valued second fundamental form $\a\in \Gamma(\Hom(TM,TM;N_fM))$. 
Let $\|\a\|^2\in C^\infty(M)$ be given at any point $x\in M^n$ by
$$
\|\a(x)\|^2=\sum_{i,j=1}^n\|\a(x)(X_i,X_j)\|^2,
$$
where $\{X_i\}_{1\leq i\leq n}$ is an orthonormal basis of $T_xM$.
Define $\phi\in C^\infty(M)$ by
\be\label{phi}
\phi^2=\frac{n}{n-1}(\|\a\|^2-n\|\mathcal{H}\|^2),
\ee
where  $\mathcal{H}$ is the mean curvature vector field of $f$. Notice that $\phi$ vanishes precisely at the umbilical points of $f$.
The metric
$$
g^*=\phi^2 g,
$$
defined on the open subset of non-umbilical points of $f$, is called the \emph{Moebius metric} determined by $f$. The metric $g^*$ is invariant under Moebius transformations of the ambient space, that is, if two immersions differ by a Moebius transformation of $\R^m$, then their corresponding Moebius metrics coincide. 

It was shown in \cite{Wa} that a hypersurface  $f\colon M^n\to\R^{n+1}$ is uniquely determined, up to Moebius transformations of the ambient space, by its Moebius metric and its \emph{Moebius shape operator} $S=\phi^{-1}(A-HI)$, where $A$ is the shape operator of $f$ with respect to a unit normal vector field $N$ and $H$ is the corresponding mean curvature function. A similar result holds for submanifolds of arbitrary codimension  (see \cite{Wa} and Section $9.8$ of \cite{DT}).

Li, Ma and Wang investigated in \cite{LMW} the natural and interesting problem of looking for the hypersurfaces $f\colon M^n\to\R^{n+1}$ that are not determined, up to Moebius transformations of $\R^{n+1}$,  only by their Moebius metrics. This fits into the fundamental problem in Submanifold theory of looking for data that are sufficient to determine a submanifold up to some group of transformations of the ambient space. 
 
  More precisely, an umbilic-free hypersurface $f\colon M^n\to\R^{n+1}$ is said to be \emph{Moebius deformable} if there exists an immersion $\tilde f\colon M^n\to\R^{n+1}$ that shares with $f$ the same Moebius metric and is not Moebius congruent to $f$ on any open subset of $M^n$. The first result in \cite{LMW} is that a Moebius deformable hypersurface with dimension $n\geq 4$ must  carry a principal curvature with multiplicity at least $n-2$.
As pointed out in \cite{LMW}, for $n\geq 5$ this is already a consequence of Cartan's classification in \cite{Ca} (see also \cite{DT1} and Chapter $17$ of~\cite{DT}) of the more general class of \emph{conformally deformable} hypersurfaces. These are the hypersurfaces  $f\colon M^n\to\R^{n+1}$ that admit a non-trivial \emph{conformal deformation} $\tilde f\colon M^n\to\R^{n+1}$, that is,  an immersion such that $f$ and $\tilde f$ induce conformal metrics on $M^n$ and do not differ by a Moebius transformation of $\R^{n+1}$ on any open subset of $M^n$ .

According to Cartan's classification, besides the conformally flat hypersurfaces, which have a principal curvature with multiplicity greater than or equal to $n-1$ and are highly conformally deformable, the remaining ones fall into one of the following classes:
\begin{itemize}
\item[(i)] \emph{conformally surface-like hypersurfaces}, that is, those that differ by a Moebius transformation of $\mathbb{R}^{n+1}$ from cylinders and rotation hypersurfaces over surfaces in $\R^3$, or from cylinders over three-dimensional 
hypersurfaces of $\R^4$ that are cones over surfaces in $\Sf^3$; 
\item[(ii)]  \emph{conformally ruled hypersurfaces}, that is, hypersurfaces $f\colon M^n\to\R^{n+1}$ for which $M^n$ carries an integrable $(n-1)$-dimensional distribution 
whose leaves are mapped by $f$ into umbilical submanifolds of $\R^{n+1}$;
\item[(iii)] hypersurfaces that admit a  non-trivial \emph{conformal variation}  $F\colon (-\epsilon,\epsilon)\times M^n \to \mathbb{R}^{n+1}$, that is,  a smooth map defined on the product of an open interval $(-\epsilon,\epsilon)\subset \mathbb{R}$ with $M^n$ such that, for any $t \in (-\epsilon,\epsilon)$, the map  $f_t = F(t; \cdot)$,  with $f_0 = f$, is a non-trivial \emph{conformal deformation} of $f$;
\item[(iv)]  hypersurfaces that admit a single non-trivial conformal deformation.
\end{itemize}

 It was shown in \cite{LMW}  that, among the {conformally surface-like hypersurfaces,  the ones that are Moebius deformable are those that  are determined by a Bonnet surface  $h\colon L^2\to \mathbb{Q}_\epsilon^3$ admitting isometric deformations preserving the mean curvature function.  Here $\mathbb{Q}_\epsilon^3$ stands for a space form of constant sectional curvature $\epsilon \in \{-1,0,1\}$. It was also shown in \cite{LMW} that an umbilic-free conformally flat hypersurface $f\colon M^n\to \R^{n+1}$, $n\geq 4$ (hence with a principal curvature of constant multiplicity $n-1$), admits non-trivial deformations preserving the Moebius metric if and only if it has constant Moebius curvature, that is, its Moebius metric has constant sectional curvature.  Such hypersurfaces were classified in \cite{GLLMW}, and an alternative proof of the classification was given in \cite{LMW}. They were shown to be, up to Moebius transformations of $\R^{n+1}$, either cylinders or rotation hypersurfaces over the so-called \emph{curvature spirals} in $\R^2$ or $\R_+^2$, respectively, the latter endowed with the hyperbolic metric, or cylinders over surfaces that are cones over curvature spirals in $\Sf^2$.
 
It is claimed in \cite{LMW} that there exists only one further example of a Moebius deformable hypersurface, which belongs to the third of the above classes in Cartan's classification of the conformally deformable hypersurfaces. Namely, the hypersurface given by
\begin{equation}\label{example}
f=\Phi\circ (\id \times f_1)\colon  M^n:=\mathbb{H}_{-m}^{n-3}\times N^3\to \R^{n+1},\,\,\,m=\sqrt{\frac{n-1}{n}},
\end{equation}
 where $\id$ is the identity map of $\mathbb{H}_{-m}^{n-3}$,  
$f_1\colon N^3\to \mathbb{S}_m^4$ is Cartan's minimal isoparametric hypersurface, which is a tube over the Veronese embedding of $\mathbb{R}\mathbb{P}^2$ into $\mathbb{S}_m^4$, and 
$\Phi\colon \mathbb{H}_{-m}^{n-3}\times \mathbb{S}_m^{4}\subset \mathbb{L}^{n-2}\times \mathbb{R}^{5}\to \mathbb{R}^{n+1}\setminus \R^{n-4}$ is the  conformal diffeomorphism given by 
$$
\Phi(x, y)=\frac{1}{x_0}(x_1, \ldots, x_{n-4},y)
$$
for all $x=x_0e_0+x_1e_1+\cdots+ x_{n-3}e_{n-3}\in \mathbb{L}^{n-2}$ and $y=(y_1, \ldots, y_{5})\in \mathbb{S}^{4}\subset \mathbb{R}^{5}$. Here  $\{e_0, \ldots, e_{n-3}\}$ denotes a pseudo-orthonormal basis of the Lorentzian space $\mathbb{L}^{n-2}$ with $\<e_0, e_0\>=0=\<e_{n-3}, e_{n-3}\>$ and $\<e_0, e_{n-3}\>=-1/2$. The deformations of $f$ preserving the Moebius metric have been shown to be  actually compositions ${f}_t=f\circ \phi_t$ of $f$ with the elements of a one-parameter family of isometries  $\phi_t\colon M^n\to M^n$  with respect to the Moebius metric; hence all of them have the same image as $f$. 

  The initial goal of this article  was to investigate the larger class of \emph{infinitesimally} Moebius bendable hypersurfaces, that is, umbilic-free hypersurfaces $f\colon M^n\to \mathbb{R}^{n+1}$ for which there exists  a one-parameter family of immersions $f_t\colon M^n\to \mathbb{R}^{n+1}$, with $t\in (-\epsilon, \epsilon)$ and $f_0=f$, such that the Moebius metrics determined by $f_t$ coincide up to the first order, in the sense that  
  $\frac{\partial}{\partial t}|_{t=0}g_t^*=0$. This is carried out for $n\geq 5$ in the forthcoming paper \cite{JT2}. 
  
  In the course of our investigation, however, we realized that the infinitesimally Moebius bendable hypersurfaces of dimension $n\geq 5$ in our classification that are not conformally surface-like are actually also Moebius deformable. Nevertheless, except for the example in the preceding paragraph, they do not appear in the classification of such hypersurfaces as stated in \cite{LMW}. This has led us to revisit that classification under a different approach from that in \cite{LMW}.

 To state our  result, we need to recall some terminology.  Let $f\colon M^n\to \R^{n+1}$ be an oriented hypersurface with respect to a unit normal vector field~$N$. 
Then the family of hyperspheres 
$
x\in M^n\mapsto S(h(x),r(x))
$
with radius $r(x)$ and center 
$
h(x)=f(x)+r(x)N(x)
$
is enveloped by $f$.  If, in particular, $1/r$ is the mean curvature of $f$, it is called the \emph{central sphere congruence} of $f$.

Let $\mathbb{V}^{n+2}$ denote the light cone in the Lorentz space $\mathbb{L}^{n+3}$ and let $\Psi=\Psi_{v,w,C}\colon\mathbb{R}^{n+1}\to\mathbb{L}^{n+3}$
be the isometric embedding onto 
$$
\mathbb{E}^{n+1}=\mathbb{E}^{n+1}_w=\{u\in\mathbb{V}^{n+2}:\<u,w\> =1\}\subset\mathbb{L}^{n+3}
$$
given by 
\be\label{eq:Psi}
\Psi(x)=v+Cx-\frac{1}{2}\|x\|^2w,
\ee
in terms of $w\in \mathbb{V}^{n+2}$, $v\in \mathbb{E}^{n+1}$ 
and a linear isometry $C\colon\mathbb{R}^{n+1}\to\{v,w\}^\perp$. Then the congruence of hyperspheres  $x\in M^n\mapsto S(h(x),r(x))$ is determined by  
the map 
$S\colon M^n\to\mathbb{S}_{1,1}^{n+2}$ that takes values in the Lorentzian sphere 
$$
\Sf_{1,1}^{n+2}=\{x\in\Les^{n+3}\colon\<x,x\>=1\}
$$
and is defined by
$$
S(x)=\frac{1}{r(x)}\Psi(h(x))+\frac{r(x)}{2}w,
$$
in the sense that $\Psi(S(h(x),r(x)))=\mathbb{E}^{n+1}\cap S(x)^\perp$ for all $x\in M^n$. 
The map $S$ has rank $0<k<n$, that is, it corresponds to a $k$-parameter congruence of hyperespheres, if and only if $\lambda=1/r$ is a principal curvature of $f$ with  constant multiplicity $n-k$ (see Section $9.3$ of \cite{DT} for details). In this case, $S$ gives rise to a map $s\colon L^k\to \mathbb{S}_{1,1}^{n+2}$ such that $S\circ \pi=s$, where $\pi\colon M^n\to L^k$ is the canonical projection onto the quotient space of leaves of $\ker (A-\lambda I)$.

  \begin{theorem}\label{thm:crux}
Let $f\colon M^n\to\R^{n+1}$, $n\geq 5$, be a Moebius deformable hypersurface  that is not conformally surface-like on  any open subset and has a principal curvature of constant multiplicity $n-2$. 
Then the central sphere congruence of $f$ is determined by a minimal space-like surface $s\colon L^2\to\Sf_{1,1}^{n+2}$.

Conversely, any simply connected hypersurface $f\colon M^n\to\R^{n+1}$, $n\geq 5$, whose central sphere congruence is determined by a minimal space-like surface $s\colon L^2\to\Sf_{1,1}^{n+2}$ is Moebius deformable. In fact, $f$ is Moebius bendable: it admits precisely a one-parameter family of conformal deformations, all of which share with $f$ the same Moebius metric. 
\end{theorem}

\begin{remarks}\label{rem}\emph{  1) Particular examples of Moebius deformable hypersurfaces  $f\colon M^n\to\R^{n+1}$ that are not conformally surface-like on  any open subset and have a principal curvature of constant multiplicity $n-2$ are the minimal hypersurfaces of rank two. These are well-known to admit a one-parameter associated family of isometric deformations, all of which are also minimal of rank two.
The elements of the associated family, sharing with $f$ the same induced metric, all have the same scalar curvature and, being minimal, also share with $f$ the same Moebius metric. These examples are not comprised in the statement of Proposition $9.2$  in \cite{LMW} and, since the elements of the associated family of a  minimal hypersurface of rank two do not have in general the same image, neither in the statement of Theorem $1.5$ therein. \vspace{1ex}\\
2) More general examples are the compositions $f=P\circ h$ of minimal hypersurfaces $h\colon M^n\to\Q_c^{n+1}$ of rank two with a ``stereographic projection" $P$ of $\Q_c^{n+1}$ (minus one point if $c> 0$) onto $\R^{n+1}$. The latter are precisely the hypersurfaces $f\colon M^n\to\R^{n+1}$  with a principal curvature of constant multiplicity $n-2$ whose central sphere congruences are determined by  minimal space-like surfaces $s\colon L^2\to\Sf_{1,1}^{n+2}\subset\mathbb{L}^{n+3}$  such that $s(L)$ is contained in a hyperplane of $\mathbb{L}^{n+3}$ orthogonal to a vector $T\in \mathbb{L}^{n+3}$ satisfying $-\<T,T\>=c$
(see, e.g., Corollary $3.4.6$ in \cite{H-J}). \vspace{1ex}\\
3) The central sphere congruence of the hypersurface given by \eqref{example} is a Veronese surface in a sphere $\mathbb{S}^4\subset \mathbb{S}_{1,1}^{n+2}$. 
\vspace{1ex}\\
4) The proof of Theorem \ref{thm:crux}  makes use of some arguments in the classification of the conformally deformable hypersurfaces of dimension $n\geq 5$  given in Chapter 17 of~\cite{DT}. 
}
\end{remarks}

\section{Preliminaries}

In this short section we recall some basic definitions and state Wang's fundamental theorem for hypersurfaces in Moebius geometry. 

Let $f\colon M\to\R^{n+1}$ be an umbilic-free immersion with Moebius metric $g^*=\<\cdot,\cdot\>^*$ and Moebius shape operator $S$. The \emph{Blaschke tensor} $\psi$ of $f$
is the endomorhism defined by 
$$
\<\psi X,Y\>^*=\frac{H}{\phi}\<SX,Y\>^*
+\frac{1}{2\phi^2}\left(\|\grad^* \phi\|_*^2
+H^2\right)\<X,Y\>^*-\frac{1}{\rho}\hess^* \phi(X,Y)
$$
for all $X,Y\in \mathfrak{X}(M)$, where $\grad^*$ and $\hess^*$ stand for the gradient and Hessian, respectively,  with respect to $g^*$.
The \emph{Moebius form}
$\omega\in\Gamma(T^*M)$ of $f$ is defined by 
$$
\omega(X)=-\frac{1}{\phi}\<\grad^*H+S\grad^*\phi,X\>^*.
$$

The Moebius shape operator, the Blaschke tensor and the Moebius form of $f$ are Moebius invariant tensors that satisfy the conformal Gauss and Codazzi equations
\be\label{gaussconf}
R^*(X,Y)=SX\wedge^*SY+\psi X\wedge^*Y+X\wedge^*\psi Y
\ee
and 
\be\label{codconf}
(\nabla^*_XS)Y-(\nabla^*_YS)X=\omega(X)Y-\omega(Y)X
\ee
for all $X,Y\in\mathfrak{X}(M)$, where  $\nabla^*$ denotes the Levi-Civita connection, $R^*$ the curvature tensor and $\wedge^*$ the wedge product with respect to $g^*$. We also point out for later use that the Moebius shape operator $S=\phi^{-1}(A-HI)$, besides being traceless, has constant norm $ \sqrt{(n-1)/n}$.

The following fundamental result was proved by Wang (see Theorem 3.1 in \cite{Wa}).

\begin{proposition}\label{congr}
Two umbilic-free hypersurfaces $f_1,f_2\colon M^n\to \R^{n+1}$ are conformally (Moebius) congruent if and only if they share the same Moebius metric and the same Moebius second fundamental form (up to sign). 
\end{proposition} 

\section{Proof of Theorem \ref{thm:crux}}

  This section is devoted to the proof of Theorem \ref{thm:crux}. In the first subsection we
   use the theory of flat bilinear forms to give an alternative proof of a key proposition proved in \cite{LMW} on the structure of the Moebius shape operators of Moebius deformable hypersurfaces. The proof of Theorem \ref{thm:crux} is provided in the subsequent subsection.

\subsection{Moebius shape operators of Moebius deformable hypersurfaces}

  The starting point for the proof of Theorem \ref{thm:crux} is Proposition \ref{kermoeb} below, which gives the structure of the Moebius shape operator of a Moebius deformable hypersurface of dimension $n\geq 5$ that carries a principal curvature of multiplicity $(n-2)$ and is not conformally surface-like on any open subset.  
  
  First we provide, for the sake of completeness, an alternative proof for $n\geq~5$, based on the theory of flat bilinear forms, of a result first proved for $n\geq 4$ by Li, Ma and Wang in \cite{LMW} (see Theorem $6.1$ therein) on the structure of the Moebius shape operators of any pair of Euclidean hypersurfaces of dimension $n\geq 5$ that are Moebius deformations of each other (see Proposition~\ref{commeig} below).  

Recall that if $W^{p,q}$ is a vector space of dimension $p+q$  endowed with an  
inner product $\<\!\<\,,\,\>\!\>$ of signature $(p,q)$, and $V$,  $U$ are
finite dimensional vector spaces, then a  bilinear form 
$\beta\colon V\times U\to W^{p,q}$ is said to be 
\emph{flat} with respect to $\<\!\<\,,\,\>\!\>$ if 
$$
\<\!\<\,\beta(X,Y),\beta(Z,T)\>\!\>-\<\!\<\,\beta(X,T),\beta(Z,Y)\>\!\>=0
$$
for all  $X,Z\in V$ and $Y,T\in U$. It is called
\emph{null} if
$$
\<\!\<\,\beta(X,Y),\beta(Z,T)\>\!\>=0
$$
for all  $X,Z\in V$ and $Y,T\in U$. Thus a null bilinear form is
necessarily flat.

\begin{proposition}\label{flat}
Let $f_1,f_2\colon M^n\to\R^{n+1}$, be  umbilic-free immersions that share the same Moebius metric $\<\, ,\,\>^*$.  Let $S_i$ and $\psi_i$, $i=1,2$, denote their corresponding Moebius shape operators and Blaschke tensors. Then, for each $x\in M^n$,  the bilinear form $\Theta\colon T_xM\times T_xM\to\R^{2,2}$ defined by 
$$
\Theta(X,Y)=(\<S_1X,Y\>^*,\frac{1}{\sqrt{2}}\<\Psi_+X,Y\>^*,\<S_2X,Y\>^*,\frac{1}{\sqrt{2}}\<\Psi_-X,Y\>^*),
$$
 where $\Psi_{\pm}=I\pm(\psi_1-\psi_2)$, is flat with respect to the (indefinite) inner product $\<\!\<\cdot,\cdot\>\!\>$ in $\R^{2,2}$. Moreover, $\Theta$  is null for all $x\in M^n$ if and only if $f_1$ and $f_2$ are Moebius congruent.
\end{proposition}
\proof Using \eqref{gaussconf}  for $f_1$ and $f_2$  we obtain
\begin{align*}
    \<\!\<\Theta(X,Y),&\Theta(Z,W)\>\!\>-\<\!\<\Theta(X,W),\Theta(Z,Y)\>\!\>\\
=&\<(S_1Z\wedge^*S_1X)Y,W\>^*-\<(S_2Z\wedge^*S_2X)Y,W\>^*\\
&+\<((\psi_1-\psi_2)Z\wedge^*X)Y,W\>^*+\<(Z\wedge^*(\psi_1-\psi_2)X)Y,W\>^*\\
=&0
\end{align*}
for all $x\in M^n$ and $X, Y, Z, W\in T_xM$, which proves the first assertion. 

Assume now that $\Theta$ is null for all $x\in M^n$. 
 Then 
\begin{align*}
 0=& \<\!\<\Theta(X,Y),\Theta(Z,W)\>\!\>= \<S_1X,Y\>^*\<S_1Z,W\>^*-\<S_2X,Y\>^*\<S_2Z,W\>^*\\
    &+\frac{1}{2}\<(I+(\psi_1-\psi_2))X,Y)\>^*
\<(I+(\psi_1-\psi_2))Z,W)\>^*\\
&-\frac{1}{2}\<(I-(\psi_1-\psi_2))X,Y)\>^*\<(I-(\psi_1-\psi_2))Z,W)\>^*
\end{align*}
for all $x\in M^n$ and $X, Y, Z, W\in T_xM$. 
This is equivalent to
\begin{align}\label{eq:s1s2}
&\;\<S_1X,Y\>^*S_1-\<S_2X,Y\>^*S_2+\frac{1}{2}\<(I+(\psi_1-\psi_2))X,Y)\>^*(I+(\psi_1-\psi_2))\nonumber\\
&-\frac{1}{2}\<(I-(\psi_1-\psi_2))X,Y)\>^*(I-(\psi_1-\psi_2))\nonumber\\
=&\<S_1X,Y\>^*S_1-\<S_2X,Y\>^*S_2+\<X,Y\>^*(\psi_1-\psi_2)+\<(\psi_1-\psi_2)X,Y\>^*I\\
=&0\nonumber
\end{align}
for all  $x\in M^n$ and $X,Y\in T_xM$. Now we use that
\begin{equation}\label{eq:prop920}
(n-2)\<\psi_i X,Y\>^*=Ric^*(X,Y)+\<S_i^2X,Y\>^*-\frac{n^2s^*+1}{2n}\<X,Y\>^*
\end{equation}
for all $X,Y\in T_xM$,  where $Ric^*$ and $s^*$ are the Ricci and scalar curvatures of the Moebius metric (see, e.g., Proposition $9.20$ in \cite{DT}), which implies that
$$
\trace \psi_1=\frac{n^2s^*+1}{2n}=\trace\psi_2.
$$
Therefore, taking traces in \eqref{eq:s1s2} yields
$$
\<(\psi_1-\psi_2)X,Y\>^*=0
$$
for all $x\in M^n$ and $X, Y\in T_xM$. Thus $\psi_1=\psi_2$, and hence
$
\<S_1X,Y\>^*S_1=\<S_2X,Y\>^*S_2.
$
In particular, $S_1$ and $S_2$ commute. Let $\lambda_i$ and $\rho_i$, $1\leq i\leq n$, denote their respective eigenvalues. Then 
$
\lambda_i\lambda_j=\rho_i\rho_j
$
for all $1\leq i,j\leq n$ and, in particular, $\lambda_i^2=\rho_i^2$ for any $1\leq i\leq n$.
If $\lambda_1=\rho_1\neq 0$, then $\lambda_j=\rho_j$ for any $j$, and then $S_1=S_2$. Similarly, if $\lambda_1=-\rho_1\neq 0$, then $S_1=-S_2$. Therefore, in any case, $f_1$ and $f_2$ are Moebius congruent by Proposition \ref{congr}.\vspace{2ex}
\qed


\begin{proposition}\label{commeig}
Let $f_1,f_2\colon M^n\to\R^{n+1}$, $n\geq5$, be umbilic-free immersions that are Moebius deformations of each other. Then there exists a distribution $\Delta$ of rank $(n-2)$ on an open and dense subset  $\mathcal{U}\subset M^n$ such that, for each $x\in \mathcal{U}$,  $\Delta(x)$ is contained in eigenspaces of the Moebius shape operators of both $f_1$ and
$f_2$ at $x$ correspondent to a common eigenvalue (up to sign).
\end{proposition}
\proof  First notice that, for each $x\in M^n$, the kernel
$$\mathcal{N}(\Theta):=\{Y\in T_xM\,:\, \Theta(X,Y)=0\,\,\mbox{for all}\,\, X\in T_xM\}$$
of the flat bilinear form $\Theta\colon T_xM\times T_xM\to\R^{2,2}$ given by Proposition \ref{flat}   is trivial, for if $Y\in T_xM$ belongs to $\mathcal{N}(\Theta)$, then $\<\Psi_+Y,Y\>^*=0=\<\Psi_-Y,Y\>^*$, which implies that $\<Y,Y\>=0$, and hence $Y=0$.

Now, by Proposition \ref{congr} and the last assertion in Proposition \ref{flat}, the flat bilinear form $\Theta$ is not null on any open subset of $M^n$, for  $f_1$ and $f_2$ are not Moebius congruent on any open subset of $M^n$. Let  $\mathcal{U}\subset M^n$ be the open and dense subset where $\Theta$ is not null.
Since $n\geq 5$, it follows from Lemma~4.22 in \cite{DT}  that, at any $x\in \mathcal{U}$, there exists an orthogonal decomposition
$
\R^{2,2}=W_1^{1,1}\oplus W_2^{1,1}
$
according to which $\Theta$ decomposes as $\Theta=\Theta_1+\Theta_2$, where $\Theta_1$ is null and $\Theta_2$ is flat with $\dim \mathcal{N}(\Theta_2)\geq n-2$. 

We claim that $\Delta=\mathcal{N}(\Theta_2)$ is contained in eigenspaces of both $S_1$ and $S_2$ at any $x\in \mathcal{U}$. In order to prove this, take any $T\in\Gamma(\Delta)$, so that $\Theta(X,T)=\Theta_1(X,T)$ for any $X\in T_xM$, and hence
$
\<\!\<\Theta(X,T),\Theta(Z,Y)\>\!\>=0
$
for all $X,Y,Z\in T_xM$. Equivalently,
\be\label{rels1s2}
\<S_1X,T\>^*S_1-\<S_2X,T\>^*S_2+\<(\psi_1-\psi_2)X,T\>^*I+\<X,T\>^*(\psi_1-\psi_2)=0
\ee
for any $X\in T_xM$. In particular, for $X$ orthogonal to $T$, 
$$
\<S_1X,T\>^*S_1-\<S_2X,T\>^*S_2+\<(\psi_1-\psi_2)X,T\>I=0.
$$
Assume that $T$ is not an eigenvector of $S_1$. Then there exists $X$ orthogonal to $T$ such that $\<S_1X,T\>^*\neq 0$. Since $f_1$ is umbilic-free, we must have $\<S_2X,T\>^*\neq 0$.
Thus $S_1$ and $S_2$ are mutually diagonalizable. Let $X_1, \ldots, X_n$ be an orthonormal diagonalizing basis of both $S_1$ and $S_2$ with respective eigenvalues $\lambda_i$ and $\rho_i$, $1\leq i\leq n$. Since $T$ is not an eigenvector, there are at least two distinct eigenvalues, say, $0\neq\lambda_1\neq\lambda_2$, with corresponding eigenvectors $X_1$ and $X_2$, such that $\<X_1,T\>^*\neq0\neq\<X_2,T\>^*$. Thus \eqref{rels1s2} yields
$$
\lambda_1\<X_1,T\>^*S_1-\rho_1\<X_1,T\>^*S_2+\<(\psi_1-\psi_2)X_1,T\>^*I+\<X_1,T\>^*(\psi_1-\psi_2)=0
$$
and
$$
\lambda_2\<X_2,T\>^*S_1-\rho_2\<X_2,T\>^*S_2+\<(\psi_1-\psi_2)X_2,T\>^*I+\<X_2,T\>^*(\psi_1-\psi_2)=0.
$$
It follows from \eqref{eq:prop920} that 
$
(n-2)(\psi_1-\psi_2)=S_1^2-S_2^2.
$
Hence
$$
\lambda_1S_1-\rho_1S_2+\frac{1}{n-2}(\lambda_1^2-\rho_1^2)I+(\psi_1-\psi_2)=0
$$
and
$$
\lambda_2S_1-\rho_2S_2+\frac{1}{n-2}(\lambda_2^2-\rho_2^2)I+(\psi_1-\psi_2)=0.
$$
Taking traces in the above expressions we obtain 
$$\lambda_1^2-\rho_1^2=0=\lambda_2^2-\rho_2^2.$$
On the other hand, the above relations also yield
$$
\lambda_1\lambda_i-\rho_1\rho_i+\frac{1}{n-2}(\lambda_i^2-\rho_i^2)=0
$$
and
$$
\lambda_2\lambda_i-\rho_2\rho_i+\frac{1}{n-2}(\lambda_i^2-\rho_i^2)=0
$$
for any $1\leq i\leq n$.
Assume first that $\lambda_1=\rho_1$, and hence $\lambda_2=\rho_2$.
Then the preceding expressions become
$$
(\lambda_i-\rho_i)\left(\lambda_j+\frac{1}{n-2}(\lambda_i+\rho_i)\right)=0
$$
for $j=1,2$ and $1\leq i\leq n$.
Since $S_1\neq S_2$ and both tensors have vanishing trace, there must exist at least two directions for which $\lambda_i-\rho_i\neq0$.
For such a fixed direction, say $k$, we have 
$$
\lambda_j+\frac{1}{n-2}(\lambda_k+\rho_k)=0,
$$
with $j=1,2$. Thus $\lambda_1=\lambda_2$, which is a contradiction.

Similarly, if we assume $\lambda_1=-\rho_1$, we obtain that $\lambda_2=-\rho_2$, and then
$$
(\lambda_i+\rho_i)\left(\lambda_j+\frac{1}{n-2}(\lambda_i-\rho_i)\right)=0
$$
for $j=1,2$ and $1\leq i\leq j$.
By the same argument as above, we see that $\lambda_1=\lambda_2$, reaching again a contradiction.
Therefore  $T$ must be an eigenvector of $S_1$. 
 Since $S_2$ is not a multiple of the identity, taking $X$ orthogonal to $T$ we see from \eqref{rels1s2} that $T$ must also be an eigenvector of $S_2$. Given that $T\in\Gamma(\Delta)$ was chosen arbitrarily, we conclude that $\Delta$ is contained in eigenspaces of both $S_1$ and $S_2$.
 
 Let $\mu_1$ and $\mu_2$ be such that $S_1|_\Delta=\mu_1 I$ and $S_2|_\Delta=\mu_2 I$. 
By \eqref{rels1s2} we have
$$
\mu_1^2-\mu_2^2+\frac{2}{n-2}(\mu_1^2-\mu_2^2)=0.
$$
Thus
$
\mu_1^2-\mu_2^2=0,
$
and hence $\mu_1=\pm \mu_2$.  

It remains to argue that $\dim\Delta=n-2$. After changing the normal vector of either $f_1$ or $f_2$, if necessary, one can assume that $\mu_1=\mu_2:=\mu$. Since $S_1|_\Delta=\mu I=S_2|_\Delta$, if  $\dim\Delta=n-1$ then the condition $\trace(S_1)=0=\trace(S_2)$ would imply that $S_1=S_2$, a contradiction. 
\vspace{1ex}\qed

Now we make the extra assumptions that $f$ is not conformally surface-like
on any open subset of $M^n$ and has a principal curvature with constant multiplicity $n-2$.
 
\begin{proposition}\label{kermoeb}
Let $f_1\colon M^n\to\R^{n+1}$, $n\geq 5$, be a Moebius deformable hypersurface with a principal curvature $\lambda$ of constant multiplicity $n-2$. Assume that $f_1$ is not conformally surface-like
on any open subset of $M^n$. If $f_2\colon M^n\to\R^{n+1}$ is a Moebius deformation of $f_1$, then the Moebius shape operators $S_1$ and $S_2$ of $f_1$ and $f_2$, respectively, have constant eigenvalues $\pm\sqrt{(n-1)/2n}$ and $0$, and the eigenspace $\Delta$ correspondent to $\lambda$ as a common kernel. In particular, $\lambda$ and the corresponding principal curvature of $f_2$ coincide with the mean curvatures of $f_1$ and $f_2$, respectively.
Moreover, the Moebius forms of $f_1$ and $f_2$ vanish on $\Delta$.
\end{proposition}

For the proof of Proposition \ref{kermoeb}, we will make use of Lemma \ref{le:split} below (see Theorem $1$ in \cite{DFT} or Corollary $9.33$ in \cite{DT}), which characterizes conformally surface-like hypersurfaces among hypersurfaces of dimension $n$ that carry a principal curvature with constant multiplicity $n-2$ in terms of the splitting tensor of the corresponding eigenbundle.
Recall that, given a distribution $\Delta$ on a Riemannian manifold $M^n$, its
\emph{splitting tensor}
$C\colon\Gamma(\Delta)\to\Gamma(\End(\Delta^\perp))$ is defined by 
$$
C_TX=-\nabla_X^hT
$$ 
for all $T\in\Gamma(\Delta)$ and $X\in\Gamma(\Delta^\perp)$, where $\nabla_X^hT=(\nabla_XT)_{\Delta^\perp}$.

\begin{lemma} \label{le:split} Let $f\colon M^n\to\R^{n+1}$, $n\geq 3$, be a hypersurface
with a principal curvature of multiplicity $n-2$ and let $\Delta$ denote its eigenbundle.
Then $f$ is  conformally surface-like if and only if the splitting tensor  of $\Delta$ satisfies $C(\Gamma(\Delta))\subset \spa\{I\}$. 
\end{lemma}

\noindent  \emph{Proof of Proposition \ref{kermoeb}:} Since  $f_1$ has a principal curvature $\lambda$ of constant multiplicity $n-2$, it follows from Proposition \ref{commeig} that, 
after changing the normal vector field of either $f_1$ or $f_2$, if necessary, we can assume that the Moebius shape operators $S_1$ and $S_2$ of $f_1$ and $f_2$ have a common eigenvalue $\mu$ with the same eigenbundle $\Delta$ of rank $n-2$.

Let $\lambda_i$, $i=1,2$, be the eigenvalues of $S_1|_{\Delta^\perp}$. In particular, $\lambda_1\neq\mu\neq\lambda_2$. The conditions $\trace(S_1)=0=\trace(S_2)$ and $\|S_1\|^2=(n-1)/n=\|S_2\|^2$ imply that $S_1$ and $S_2$ have the same eigenvalues. Then we must also have $\lambda_1\neq\lambda_2$, for otherwise $S_1$ and $S_2$ would coincide.

Let $X,Y\in\Gamma(\Delta^\perp)$ be an orthonormal frame of eigenvectors of $S_1|_{\Delta^\perp}$ with respect to $g^*$. Then $S_1X=\lambda_1X$, $S_1Y=\lambda_2Y$,  $S_2X=b_1X+cY$ and $S_2Y=cX+b_2Y$ for some smooth functions $b_1$, $b_2$ and $c$.
Since $\trace(S_1)=0=\trace(S_2)$ and $\|S_1\|^{*2}=(n-1)/n=\|S_2\|^{*2}$, we have
\begin{align}
\lambda_1+\lambda_2+(n-2)\mu&=0\label{trace1},\\
\lambda_1^2+\lambda_2^2+(n-2)\mu^2&=\frac{n-1}{n},\label{norm1}\\
b_1+b_2+(n-2)\mu&=0\label{trace2},\\
b_1^2+b_2^2+2c^2+(n-2)\mu^2&=\frac{n-1}{n}.\label{norm2}
\end{align}

  Thus the first assertion in the statement will be proved once we show that $\mu$ vanishes identically. The last assertion will then be an immediate consequence of \eqref{codconf}.

The umbilicity of $\Delta$, together with \eqref{codconf} evaluated in orthonormal sections $T$ and $S$ of $\Delta$ with respect to $g^*$, imply that
$\omega_1(T)=T(\mu)=\omega_2(T)$, where $\omega_i$ is the Moebius form of $f_i$, $1\leq i\leq 2$. Taking the derivative of \eqref{trace1} and \eqref{norm1} with respect to $T\in\Gamma(\Delta)$,  we obtain
$$
T(\lambda_1)=\frac{(n-2)(\mu-\lambda_2)}{\lambda_2-\lambda_1}T(\mu)\;\;\mbox{and}\;\;
T(\lambda_2)=\frac{(n-2)(\lambda_1-\mu)}{\lambda_2-\lambda_1}T(\mu).
$$

The $X$ and $Y$ components of \eqref{codconf} for $S_1$ evaluated in $X$ and $T\in\Gamma(\Delta)$ give, respectively,
\be\label{cod1xtx}
(\mu-\lambda_1)\<\nabla^*_XT,X\>^*=T(\lambda_1)-T(\mu)=-\frac{n\lambda_2}{\lambda_2-\lambda_1}T(\mu)
\ee
and
\be\label{cod1xty}
(\mu-\lambda_2)\<\nabla^*_XT,Y\>^*=(\lambda_1-\lambda_2)\<\nabla^*_TX,Y\>^*.
\ee
Similary, the $X$ and $Y$ components of \eqref{codconf} for $S_1$ evaluated in $Y$ and $T$ give, respectively,
\be\label{cod1ytx}
(\mu-\lambda_1)\<\nabla^*_YT,X\>^*=(\lambda_2-\lambda_1)\<\nabla^*_TY,X\>^*
\ee
and
\be\label{cod1yty}
(\mu-\lambda_2)\<\nabla^*_YT,Y\>^*=T(\lambda_2)-T(\mu)=\frac{n\lambda_1}{\lambda_2-\lambda_1}T(\mu).
\ee

We claim that $S_1$ and $S_2$ do not commute, that is, that $c\neq 0$. Assume otherwise. Then Eqs. \eqref{trace1} to \eqref{norm2} imply that $S_2X=\lambda_2X$ and $S_2Y=\lambda_1Y$. Hence, the $X$ and $Y$ components of \eqref{codconf} for $S_2$ evaluated in $X$ and $T\in\Gamma(\Delta)$ give, respectively,
\be\label{cod2xtx}
(\mu-\lambda_2)\<\nabla^*_XT,X\>^*=T(\lambda_2)-T(\mu)
\ee
and
\be\label{cod2xty}
(\mu-\lambda_1)\<\nabla^*_XT,Y\>^*=(\lambda_2-\lambda_1)\<\nabla^*_TX,Y\>^*.
\ee
Similary, the $X$ and $Y$ components of \eqref{codconf} for $S_2$ evaluated in $Y$ and $T$ give, respectively,
\be\label{cod2ytx}
(\mu-\lambda_2)\<\nabla^*_YT,X\>^*=(\lambda_1-\lambda_2)\<\nabla^*_TY,X\>^*
\ee
and
\be\label{cod2yty}
(\mu-\lambda_1)\<\nabla^*_YT,Y\>^*=T(\lambda_1)-T(\mu).
\ee
Adding \eqref{cod1xty} and \eqref{cod2xty} yields
$$
(2\mu-\lambda_1-\lambda_2)\<\nabla^*_XT,Y\>^*=0.
$$
Similarly, Eqs \eqref{cod1ytx} and \eqref{cod2ytx} give
$$
(2\mu-\lambda_1-\lambda_2)\<\nabla^*_YT,X\>^*=0.
$$
If $(2\mu-\lambda_1-\lambda_2)$ does not vanish identically,  there exists an open subset $U\subset M^n$ where $\<\nabla^*_XT,Y\>^*=0=\<\nabla^*_YT,X\>^*$.
Now, from \eqref{cod1xtx} and \eqref{cod2xtx} we obtain
$$
(\lambda_2-\lambda_1)\<\nabla^*_XT,X\>^*=T(\lambda_1-\lambda_2).
$$
Similarly, using \eqref{cod1yty} and \eqref{cod2yty} we have
$$
(\lambda_1-\lambda_2)\<\nabla^*_YT,Y\>^*=T(\lambda_2-\lambda_1).
$$
The preceding equations imply that the splitting tensor $C^*$ of $\Delta$  with respect to the Moebius metric satisfies $C^*_T\in\spa\{I\}$ for any $T\in\Gamma(\Delta|_U)$. From the relation between the Levi-Civita connections of conformal metrics we obtain
\be\label{confmetrics}
C^*_T=C_T -T(\log\phi)\,I,
\ee
where $\phi$ is the conformal factor of $g^*$ with respect to the metric induced by $f_1$ and $C$ is the splitting tensor of $\Delta$ corresponding to the latter metric. Therefore, we also have $C_T\in\spa\{I\}$ for any $T\in\Gamma(\Delta|_U)$, and hence $f_1|U$ is conformally surface-like  by Lemma \ref{le:split}, a contradiction.  Thus  $(2\mu-\lambda_1-\lambda_2)$ must vanish everywhere, which, together with \eqref{trace1}, implies that also $\mu$ is everywhere vanishing. Hence $\lambda_1=-\lambda_2$, and therefore $S_1=-S_2$, which is again a contradiction, and proves the claim.

Now we compute
\begin{align*}
    \<(\nabla^*_TS_2)X,X\>^*&=\<\nabla^*_T(b_1X+cY),X\>^*-\<S_2\nabla^*_TX,X\>^*\\
    &=T(b_1)+c\<\nabla^*_TY,X\>^*-c\<\nabla^*_TX,Y\>^*\\
    &=T(b_1)+2c\<\nabla^*_TY,X\>^*.
\end{align*}
In a similar way,
\begin{align*}
    \<(\nabla^*_TS_2)Y,Y\>^*&=T(b_2)+2c\<\nabla^*_TX,Y\>^*.
\end{align*}
Adding the preceding equations and using \eqref{trace2} yield
\be\label{derttr2}
\<(\nabla^*_TS_2)X,X\>^*+\<(\nabla^*_TS_2)Y,Y\>^*=(2-n)T(\mu).
\ee
From \eqref{codconf} we obtain
\begin{align*}
    \<(\nabla^*_TS_2)X,X\>^*&=\<(\nabla^*_XS_2)T,X\>^*+T(\mu)\\
    &=\mu\<\nabla^*_XT,X\>^*-\<\nabla^*_XT,S_2X\>^*+T(\mu)\\
    &=(\mu-b_1)\<\nabla^*_XT,X\>^*-c\<\nabla^*_XT,Y\>^*+T(\mu),
\end{align*}
and similarly,
$$
\<(\nabla^*_TS_2)Y,Y\>^*=(\mu-b_2)\<\nabla^*_YT,Y\>^*-c\<\nabla^*_YT,X\>^*+T(\mu).
$$
Substituting the preceding expressions in \eqref{derttr2} gives
\be\label{derttr22}
    nT(\mu)+(\mu-b_1)\<\nabla^*_XT,X\>^*+(\mu-b_2)\<\nabla^*_YT,Y\>^*
    =c\<\nabla^*_XT,Y\>^*+c\<\nabla^*_YT,X\>^*.
\ee
Let us first focus on the terms on the left-hand side of the above equation. Using \eqref{cod1xtx} and \eqref{cod1yty} we obtain
\begin{align*}
   nT(\mu)&+(\mu-b_1)\<\nabla^*_XT,X\>^*+(\mu-b_2)\<\nabla^*_YT,Y\>^*\\
   &=nT(\mu)-\frac{n\lambda_2(\mu-b_1)}{(\mu-\lambda_1)(\lambda_2-\lambda_1)}T(\mu)
   +\frac{n\lambda_1(\mu-b_2)}{(\mu-\lambda_2)(\lambda_2-\lambda_1)}T(\mu)\\
&=\frac{(n-1)(\lambda_1-b_1)}{(\mu-\lambda_2)(\lambda_2-\lambda_1)}T(\mu).   
\end{align*}
For the right-hand side of \eqref{derttr22}, using \eqref{cod1xty} and \eqref{cod1ytx} we have
\begin{align*}
    c(\<\nabla^*_XT,Y\>^*&+\<\nabla^*_YT,X\>^*)=c\left(\frac{\lambda_1-\lambda_2}{\mu-\lambda_2}\<\nabla^*_TX,Y\>^*+\frac{\lambda_2-\lambda_1}{\mu-\lambda_1}\<\nabla^*_TY,X\>^*\right)\\
    =&c\frac{(\lambda_1-\lambda_2)(\mu-\lambda_1+\mu-\lambda_2)}{(\mu-\lambda_1)(\mu-\lambda_2)}\<\nabla^*_TX,T\>^*\\
    =&c\frac{n\mu(\lambda_1-\lambda_2)}{(\mu-\lambda_1)(\mu-\lambda_2)}\<\nabla^*_TX,Y\>^*.
\end{align*}
Therefore \eqref{derttr2} becomes
\be\label{A}
(n-1)(b_1-\lambda_1)T(\mu)=nc\mu(\lambda_1-\lambda_2)^2\<\nabla^*_TX,T\>^*.
\ee 
Now evaluate \eqref{codconf} for $S_2$ in $X$ and $T$. More specifically, the $Y$ component of that equation is
$$
    T(c)=(\mu-b_2)\<\nabla^*_XT,Y\>^*-c\<\nabla^*_XT,X\>+(b_2-b_1)\<\nabla^*_TX,Y\>^*.
$$
Substituting \eqref{cod1xtx} and \eqref{cod1xty} in the above equation, and using \eqref{trace1} and \eqref{trace2}, we obtain
\begin{align}\label{tc1}
    T(c)=&\frac{(\mu-b_2)(\lambda_1-\lambda_2)}{\mu-\lambda_2}\<\nabla^*_TX,Y\>^*+\frac{cn\lambda_2}{(\mu-\lambda_1)(\lambda_2-\lambda_1)}T(\mu)\nonumber\\
    &+(b_2-b_1)\<\nabla^*_TX,Y\>^*\nonumber\\
    =&\frac{\mu\lambda_1-\mu\lambda_2-b_2\lambda_1+b_2\lambda_2+\mu b_2-\mu b_1-\lambda_2b_2+b_1\lambda_2}{\mu-\lambda_2}\<\nabla^*_TX,Y\>^*\nonumber\\
    &+\frac{cn\lambda_2}{(\mu-\lambda_1)(\lambda_2-\lambda_1)}T(\mu)\nonumber\\
    =&\frac{n\mu(\lambda_1-b_1)}{\mu-\lambda_2}\<\nabla^*_TX,Y\>^*
    +\frac{cn\lambda_2}{(\mu-\lambda_1)(\lambda_2-\lambda_1)}T(\mu).
\end{align}
Similarly, the $X$ component of \eqref{codconf} for $S_2$ evaluated in $Y$ and $T$ gives
$$
T(c)=(\mu-b_1)\<\nabla^*_YT,X\>^*-c\<\nabla^*_YT,Y\>^*+(b_2-b_1)\<\nabla^*_TX,Y\>^*.
$$
Substituting \eqref{cod1ytx} and \eqref{cod1yty} in the above equation we obtain
\be\label{tc2}
T(c)=-\frac{cn\lambda_1}{(\mu-\lambda_2)(\lambda_2-\lambda_1)}T(\mu)+\frac{n\mu(\lambda_1-b_1)}{\mu-\lambda_1}\<\nabla^*_TX,Y\>^*.
\ee
Using \eqref{norm1}, it follows from  \eqref{tc1} and \eqref{tc2} that
\be\label{B}
(n-1)cT(\mu)=n\mu(\lambda_1-b_1)(\lambda_1-\lambda_2)^2\<\nabla^*_TX,Y\>^*.
\ee
Comparing \eqref{A} and \eqref{B} yields
$$
\mu((\lambda_1-b_1)^2+c^2)\<\nabla^*_TX,Y\>^*=0.
$$
Since $(\lambda_1-b_1)^2+c^2\neq 0$, for otherwise the immersions would be Moebius congruent, then
$
\mu\<\nabla^*_TX,Y\>^*=0.
$

If  $\mu$ does not vanish identically, then there is an open subset $U$ where $\<\nabla^*_TX,Y\>^*=0$ for any $T\in\Gamma(\Delta)$. Then \eqref{cod1xty} and \eqref{cod1ytx} imply that the splitting tensor of $\Delta$ with respect to the Moebius metric satisfies $C^*_T\in\spa\{I\}$ for any $T\in\Gamma(\Delta)$. As before, this implies that the splitting tensor of $\Delta$ with respect to the metric induced by $f_1$ also satisfies $C_T\in\spa\{I\}$ for any $T\in\Gamma(\Delta)$, and hence $f_1|_U$ is conformally surface-like  by Lemma \ref{le:split}, a contradiction. Thus  $\mu$ must vanish identically. 
\qed

\subsection{Proof of Theorem \ref{thm:crux}}

  In this subsection we prove Theorem \ref{thm:crux}. First we recall one further definition.

  Let  $f\colon M^n\to\R^{n+1}$, $n\geq 3$, be a hypersurface that carries a  
principal curvature of constant multiplicity $n-2$ with 
corresponding eigenbundle $\Delta$. Let $C\colon\Gamma(\Delta)\to\Gamma(\End(\Delta^\perp))$ be the splitting tensor of $\Delta$. Then $f$    
is said to be \emph{hyperbolic} (respectively,
\emph{parabolic} or \emph{elliptic}) if there exists 
$J\in\Gamma(\End(\Delta^\perp))$ satisfying the 
following conditions:
\begin{itemize}
\item[(i)] $J^2=I$ and $J\neq I$ (respectively, $J^2=0$, with $J\neq 0$, 
or $J^2=-I$),
\item[(ii)] $\nabla^h_T J=0$ for all $T\in\Gamma(\Delta)$,
\item[(iii)] $C(\Gamma(\Delta))\subset \spa\{I,J\}$, but   $C(\Gamma(\Delta))\not\subset \spa\{I\}$.
\end{itemize}

\noindent \emph{Proof of Theorem \ref{thm:crux}:} Let $f_2\colon M^n\to\R^{n+1}$ be a Moebius deformation of $f_1:=f$. By Proposition \ref{kermoeb},  the Moebius shape operators  $S_1$ and $S_2$ of $f_1$ and $f_2$, respectively, share a common kernel  $\Delta$ of dimension $n-2$. Let $S_i$, $i=1,2$, denote also the restriction $S_i|_{\Delta^\perp}$ and define $D\in\Gamma(\End(\Delta^\perp))$ by
$$
D=S_1^{-1}S_2.
$$
It follows from Proposition \ref{kermoeb} that  $\det D=1$ at any point of $M^n$, while Proposition \ref{congr} implies that $D$ cannot be  the identity endomorphism up to sign on any open subset $U\subset M^n$, for  otherwise $f_1|_U$ and $f_2|_U$ would be Moebius congruent by Lemma \ref{le:split}.
Therefore, we can write $D=aI+bJ$,  where $b$ does not vanish on any open subset of $M^n$ and $J\in\Gamma(\End(\Delta^\perp))$  satisfies $J^2=\epsilon I$, with $\epsilon\in \{1,0,-1\}$, $J\neq I$ if $\epsilon=1$ and $J\neq 0$ if $\epsilon=0$.

From the symmetry of $S_2$ and the fact that $b$ does not vanish on any open subset of $M^n$, we see that  $S_1J$ must be symmetric. Moreover, given that $\trace S_1=0=\trace S_2$, also $\trace S_1J=0$.

Assume first that $J^2=0$. Let $X,Y\in\Gamma(\Delta^\perp)$ be orthogonal vector fields, with $Y$ of unit length (with respect to the Moebius metric $g^*$), such that $JX=Y$ and $JY=0$. Replacing $J$ by  $\|X\|^*J$, if necessary,  we can assume that also $X$ has unit length.
Let $\alpha, \beta, \gamma\in C^{\infty}(M)$ be such that $S_1X=\alpha X+\beta Y$ and $S_1Y=\beta X+\gamma Y$, so that $S_1JX=\beta X+\gamma Y$ and $S_1JY=0$. 
From the symmetry of $S_1J$ and the fact that $\trace S_1J=0$ we obtain $\beta=0=\gamma$, and hence $\alpha=\trace S_1=0$. Thus $S_1=0$, which is a contradiction.

Now assume that $J^2=I$, $J\neq I$. Let $X,Y$ be a frame of unit vector fields (with respect to $g^*$) satisfying $JX=X$ and $JY=-Y$. Write $S_1X=\alpha X+\beta Y$ and $S_1 Y=\gamma X+\delta Y$ for some $\alpha, \beta, \gamma, \delta\in C^{\infty}(M)$, so that $S_1JX=\alpha X+\beta Y$ and $S_1JY=-\gamma X-\delta Y$.
Since  $\trace S_1J=0=\trace S_1$, then $\alpha=0=\delta$. The symmetry of $S_1$ and $S_2J$ implies that $\beta=0=\gamma$, which is again a contradiction.

Therefore, the only possible case is that $J^2=-I$. Let $X,Y\in\Gamma(\Delta^\perp)$ be a frame of unit vector fields such that $JX=Y$ and $JY=-X$.  Write as before $S_1X=\alpha X+\beta Y$ and $S_1 Y=\gamma X+\delta Y$ for some $\alpha, \beta, \gamma, \delta\in C^{\infty}(M)$. Then 
 $S_1JX=\gamma X+\delta Y$ and $S_1JY=-\alpha X-\beta Y$, hence
 $\beta=\gamma$, for $\trace S_1J=0$. From the symmetry of $S_1$ we obtain
$$
\<S_1JX,Y\>=\<JX,S_1Y\>=\<Y,S_1Y\>=\gamma\<X,Y\>+\delta=\beta\<X,Y\>+\delta, 
$$
and similarly,
$$
\<S_1JY,X\>=-\alpha-\beta\<X,Y\>.
$$
Comparing the two preceding equations, taking into account the symmetry of $S_1J$ and the fact that $\trace S_1=0$, yields
$
\beta\<X,Y\>=0.
$
If $\beta$ is nonzero, then $X$ and $Y$ are orthogonal to each other. This is also the case  if $\beta$, hence also $\gamma$, is zero, for in this case $X$ and $Y$ are eigenvectors of $S_1$. Thus, in any case, we conclude that $J$ acts as a rotation of angle $\pi/2$ on $\Delta^\perp$.

Eq. \eqref{codconf} and the fact that $\omega_i|_\Delta=0$ imply that 
the splitting tensor of $\Delta$ with respect to the Moebius metric satisfies
$$
\nabla_T^{*h}S_i=S_iC^*_T
$$
for all  $T\in\Gamma(\Delta)$ and $1\leq i\leq 2$, where $$(\nabla_T^{*h}S_i)X=\nabla_T^{*h}S_iX-S_i\nabla_T^{*h}X$$ for all $X\in\Gamma(\Delta^\perp)$ and $T\in\Gamma(\Delta)$. Here
$
\nabla^{*h}_T X=(\nabla^*_TX)_{\Delta^\perp}.
$
In particular, 
$$S_iC^*_T={C^*_T}^tS_i, \,\,1\leq i\leq 2.$$
 Therefore
$$
S_1DC^*_T=S_2C^*_T={C^*_T}^tS_2={C^*_T}^tS_1D=S_1C^*_TD,
$$
and hence 
$$
[D,C^*_T]=0.
$$
 This implies that $C^*_T$ commutes with $J$, and hence $C^*_T\in\spa\{I,J\}$ for any $T\in\Gamma(\Delta)$. It follows from \eqref{confmetrics} that also the splitting tensor $C$  of $\Delta$ corresponding to the metric induced by $f$ satisfies $C_T\in\spa\{I,J\}$ for any $T\in\Gamma(\Delta)$. Moreover, by  Lemma \ref{le:split} and the assumption that $f$ is not surface-like on any open subset,  we see that $C(\Gamma(\Delta))\not\subset \spa\{I\}$ on any open subset. Now, since $J$ acts as a rotation of angle $\pi/2$ on $\Delta^\perp$, then $\nabla_T^hJ=0$. We conclude that $f$ is  elliptic with respect to $J$. 

  By Proposition \ref{kermoeb}, the central sphere congruence $S\colon M^n\to \Sf_{1,1}^{n+2}$ of $f$ is a two-parameter congruence of hyperspheres, which therefore gives rise to a surface $s\colon L^2\to\Sf_{1,1}^{n+2}$ such that
$s=S\circ \pi$, where  $\pi\colon M^n\to L^2$ is the (local) quotient map onto the space of leaves of $\Delta$.  Since $\nabla_T^hJ=0=[C_T,J]$ for any $T\in\Gamma(\Delta)$, it follows from Corollary $11.7$ in \cite{DT} that $J$ is projectable with respect to $\pi$, that is, there exists $\bar{J} \in \mbox{End}(TL)$ such that $\bar{J}\circ \pi_*=\pi_*\circ J$. In particular, the fact that $J^2=-I$ implies that $\bar{J}^2=-I$, where we denote also by $I$ the identity endomorphism of $TL$.

  Now observe that, since $f_2$ shares with $f_1$ the same Moebius metric, its induced metric is conformal to the metric induced by $f_1$. Moreover,  $f_2$ is not Moebius congruent to $f_1$ on any open subset of $M^n$ and $f_1$ has a principal curvature of constant multiplicity $(n-2)$. Thus $f_1$ is a so-called \emph{Cartan hypersurface}.   By the proof of the classification of Cartan hypersurfaces given in  Chapter~17 of \cite{DT} (see Lemma $17.4$ therein),  the surface $s$ is \emph{elliptic} with respect to $\bar{J}$, that is, for all $\bar{X},\bar{Y}\in \mathfrak{X}(L)$ we have  
\be\label{sffsur}
\alpha^s(\bar{J}\bar{X},\bar{Y})=\alpha^s(\bar{X},\bar{J}\bar{Y}).
\ee

 We claim that $\bar{J}$ is an orthogonal tensor, that is, it acts as a rotation of angle $\pi/2$ on each tangent space of $L^2$. The minimality of $s$ will then follow from this, the fact that $\bar{J}^2=-I$ and \eqref{sffsur}.

  In order to show the orthogonality of $\bar{J}$, we use the fact that the metric $\<\cdot,\cdot\>'$ on $L^2$ induced by $s$ is related to the metric of $M^n$ by
\be\label{relmetr}
\<\bar{Z},\bar{W}\>'=\<(A-\lambda I)Z, (A-\lambda I)W\>
\ee
for all $\bar{Z},\bar{W}\in\mathfrak{X}(L)$, where $A$ is the shape operator of $f$, $\lambda$ is the principal curvature of $f$ having $\Delta$ as its eigenbundle, which coincides with the mean curvature $H$ of $f$ by Proposition~\ref{kermoeb}, and $Z$, $W$ are the horizontal lifts of $\bar{Z}$ and $\bar{W}$, respectively. Notice that  $(A-\lambda I)$ is a multiple of $S_1$. Since $S_1J$ is symmetric, then also $(A-\lambda I)J$ is symmetric. Therefore, given any $\bar{X}\in\mathfrak{X}(L)$ and denoting by $X\in\Gamma(\Delta^\perp)$ its horizontal lift, we have
\begin{align*}
    \<\bar{X},\bar{J}\bar{X}\>'&=\<(A-\lambda I)X, (A-\lambda I)JX\>\\
    &=\<(A-\lambda I)J(A-\lambda I)X,X\>\\
    &=\<J(A-\lambda I)X,(A-\lambda I)X\>\\
    &=0,
\end{align*}
where in the last step we have used that $J$ acts as a rotation of angle $\pi/2$ on $\Delta^\perp$. Using again the symmetry of $(A-\lambda I)J$, the proof of the orthogonality of $\bar{J}$ is completed by noticing that
\begin{align*}
    \<\bar{J}\bar{X},\bar{J}\bar{X}\>'&=\<(A-\lambda I)JX,(A-\lambda I)JX\>\\
    &=\<J(A-\lambda I)JX,(A-\lambda I)X\>\\
    &=\<JJ^t(A-\lambda I)X,(A-\lambda I)X\>\\
    &=-\<J^2(A-\lambda I)X,(A-\lambda I)X\>\\
    &=\<\bar{X},\bar{X}\>'.
\end{align*}

    Conversely, assume that the central sphere congruence of  $f\colon M^n\to\R^{n+1}$, with $M^n$ simply connected, is  determined by a space-like minimal surface $s\colon L^2\to \Sf_{1,1}^{n+2}$. Let $\bar{J}\in \Gamma(\End(TL))$ represent a rotation of angle $\pi/2$ on each tangent space. Then ${\bar J}^2=-I$ and the second fundamental form of $s$ satisfies \eqref{sffsur} by the minimality of $s$. In particular,  $s$ is elliptic with respect to $\bar{J}$.
By Lemma $17.4$ in \cite{DT}, the hypersurface $f$ is elliptic with respect to the lift $J\in\Gamma(\End (\Delta^\perp))$ of $\bar{J}$, where $\Delta$ is the eigenbundle correspondent to the principal curvature $\lambda$ of $f$ with multiplicity $n-2$, which coincides with its mean curvature. Therefore, the splitting tensor of $\Delta$ satisfies $C_T\in\spa\{I,J\}$ for any $T\in\Gamma(\Delta)$. Since $(A-\lambda I)C_T$ is symmetric for any $T\in \Gamma(\Delta)$, as follows from the Codazzi equation, and $C(\Gamma(\Delta))\not \subset \spa\{I\}$ on any open subset, for $f$ is not conformally surface-like on any open subset, then $(A-\lambda I)J$ is also symmetric. 

  By Theorem $17.5$ in \cite{DT}, the set of conformal deformations of $f$ is in one-to-one correspondence with the set of tensors $\bar{D}\in \Gamma(\End(TL))$ with $\det \bar{D}=1$ that satisfy the Codazzi equation 
  $$\left(\nabla'_{\bar{X}}\bar{D}\right)\bar{Y}
-\left(\nabla'_{\bar{Y}}\bar{D}\right){\bar{X}}=0$$ for all 
$\bar{X},\bar{Y}\in\mathfrak{X}(L)$, 
where $\nabla'$ is the Levi-Civita connection  of the metric induced by  $s$.
For a general elliptic hypersurface, this set either consists of a one-parameter family   (continuous class) or of a single element (discrete class; see Section $11.2$ and Exercise $11.3$ in \cite{DT}). 
The surface $s$ is then said to be of the complex type of first or second species, respectively.
For a  minimal surface $s\colon L^2\to \Sf_{1,1}^{n+2}$, each tensor $\bar{J}_\theta=\cos\theta I+\sin \theta \bar{J}$, $\theta\in [0,2\pi)$, satisfies both the condition  $\det \bar{J}_\theta=1$ and the Codazzi equation, since it is a parallel tensor in $L^2$. Thus
$\{\bar{J}_\theta\}_{\theta\in [0,2\pi)}$ is \emph{the} one-parameter family of tensors in $L^2$ having determinant one and satisfying the Codazzi equation. In particular, the surface $s$ is  of the complex type of first species.  Therefore,  the hypersurface $f$ admits a one-parameter family of conformal deformations, each of which determined by one of the tensors  $\bar{J}_\theta\in \End(TL)$, $\theta\in [0,2\pi)$.  The proof of Theorem \eqref{thm:crux} will be completed once we prove that any of such conformal deformations shares with $f$ the same Moebius metric.

  Let  $f_\theta\colon M^n\to\R^{n+1}$ be the conformal deformation of $f$ determined by $\bar{J}_\theta$. Let $F_\theta\colon M^n\to \mathbb{V}^{n+2}$ be the \emph{isometric light-cone representative} of $f_\theta$, that is, $F_\theta$ is the isometric immersion of $M^n$ into the light-cone $\mathbb{V}^{n+2}\subset \mathbb{L}^{n+3}$ given by $F_\theta=\varphi_\theta^{-1}(\Psi\circ f_\theta)$, where $\varphi_\theta$ is the conformal factor of the metric $\<\cdot,\cdot\>_\theta$ induced by $f_\theta$ with respect to the metric $\<\cdot,\cdot\>$ of $M^n$, that is, $\<\cdot,\cdot\>_\theta=\varphi_\theta^2\<\cdot,\cdot\>$, and $\Psi\colon \mathbb{R}^n\to \mathbb{V}^{n+2}$ is the isometric embedding of $\mathbb{R}^n$ into $\mathbb{V}^{n+2}$ given by~\eqref{eq:Psi}.
As shown in the proof of Lemma $17.2$ in \cite{DT}, as part of the proof of the classification of Cartan hypersurfaces of dimension $n\geq 5$ given in Chapter~17 therein, the second fundamental form of $F_\theta$ is given by
\be\label{sffF}
\a^{F_\theta}(X,Y)=\<AX,Y\>\mu-\<(A-\lambda I)X,Y\>\zeta+\<(A-\lambda I)J_\theta X,Y\>\bar{\zeta}
\ee
for all $X,Y\in\mathfrak{X}(M)$, where $\{\mu,\zeta,\bar{\zeta}\}$ is an orthonormal frame of the normal bundle of $F_\theta$ in $\Les^{n+3}$ with $\mu$ space-like, $\lambda=-\<\mu, F_\theta\>^{-1}$ and $\zeta=\lambda F_\theta+\mu$ (hence $\<\zeta,\zeta\>=-1$). Here  $J_\theta$ is the horizontal lift of $\bar{J}_\theta$, which has been extended to $TM$ by setting ${J_\theta}|_\Delta=I$.

Let $\bar{X},\bar{Y}\in\mathfrak{X}(L)$ be an orthonormal frame  such that $\bar{J}\bar{X}=\bar{Y}$ and $\bar{J}\bar{Y}=-\bar{X}$, and let $X,Y\in\Gamma(\Delta^\perp)$ be the respective horizontal lifts. It follows from \eqref{relmetr} that $\{(A-\lambda I)X,(A-\lambda I)Y\}$ is an orthonormal frame of $\Delta^\perp$.  From the symmetry of $(A-\lambda I)J$ and $(A-\lambda I)$ we have
\begin{align*}
\<J(A-\lambda I)X,(A-\lambda I)X\>&=\<(A-\lambda I)J(A-\lambda I)X,X\>\\
&=\<(A-\lambda I)X,(A-\lambda I)J X\>\\
&=\<\bar{X},\bar{J}\bar{X}\>'\\
&=0.
\end{align*}
In a similar way one verifies that $\<J(A-\lambda I)Y,(A-\lambda I)Y\>=0$ and
$$\<J(A-\lambda I)Y,(A-\lambda I)X\>=1=-\<J(A-\lambda I)X,(A-\lambda I)Y\>.$$
 Thus $J$ acts on $\Delta^\perp$ as a rotation of angle $\pi/2$. The symmetry of both $(A-\lambda I)J$ and $(A-\lambda I)$ implies that $\trace(A-\lambda I)=0=\trace (A-\lambda I)J$, hence
\be\label{traces}
\trace (A-\lambda I)J_\theta=0
\ee
for all $\theta\in [0,2\pi)$. 

  Now we use the relation between the second fundamental forms of $f_\theta$ and $F_\theta$, given  by Eq. $9.32$ in \cite{DT}, namely,
\be\label{sffF2}
\a^{F_\theta}(X,Y)=\<\varphi (A_\theta-H_\theta I)X,Y\>_2\tilde{N}-\psi(X,Y)F_\theta-\<X,Y\>\zeta_2,
\ee
where $\<\cdot,\cdot\>_\theta=\varphi_\theta^2\<\cdot,\cdot\>$ is the metric induced by $f_\theta$, $A_\theta$ and $H_\theta$ are its shape operator and mean curvature, respectively, $\psi$ is a certain symmetric bilinear form, 
 $\tilde{N}\in \Gamma(N_FM)$, with $\<\tilde{N}, F_\theta\>=0$,  is a unit space-like vector field, 
 and  $\zeta_2\in \Gamma(N_FM)$ satisfies $\<\tilde{N},\zeta_2\>=0=\<\zeta_2,\zeta_2\>$ and $\<F_\theta,\zeta_2\>=1$. 
Eqs. \eqref{sffF} and \eqref{sffF2} give 
\begin{align*}
\<(A-\lambda I)J_\theta X,Y\>&=\<\a^{F_\theta}(X,Y),\bar{\zeta}\>\\
&=\varphi_\theta\<(A_\theta-H_\theta I)X,Y\>\<\tilde{N},\bar{\zeta}\>-\<X,Y\>\<\zeta_2,\bar{\zeta}\>
\end{align*}
for all $X,Y\in\mathfrak{X}(M)$, or equivalently,
\begin{equation}\label{eq:imp}
(A-\lambda I)J_\theta=\varphi_\theta\<\tilde{N},\bar{\zeta}\>(A_\theta-H_\theta I)-\<\zeta_2,\bar{\zeta}\>I.
\end{equation}
Using that 
$$\trace (A-\lambda I)J_\theta=0=\trace (A_\theta- H_\theta I),$$
 we obtain from the preceding equation that $\<\zeta_2,\bar{\zeta}\>=0$. Thus $\bar{\zeta}\in \spa\{F_\theta,\zeta_2\}^\perp$, and hence $\bar{\zeta}=\pm\tilde{N}$. Therefore, \eqref{eq:imp} reduces to
\be\label{finalkey}
(A-\lambda I)J_\theta=\pm\varphi(A_\theta-H_\theta I).
\ee
In particular, $(A_\theta-H_\theta I)|_\Delta=0$, hence also $S_\theta|_\Delta=0$, where  $S_\theta=\phi_\theta^{-1}(A_\theta-H_\theta I)$ is the Moebius shape operator of $f_\theta$, with $\phi_\theta$  given by \eqref{phi} for $f_\theta$. Since the Moebius shape operator of an umbilic-free immersion is traceless and has constant  norm  $\sqrt{(n-1)/n}$, then $S_\theta$ must have constant eigenvalues $\sqrt{(n-1)/2n}$, $-\sqrt{(n-1)/2n}$ and $0$. The same holds for the Moebius second fundamental form $S_1$ of $f$, which has also $\Delta$ as its kernel. We conclude that the eigenvalues of $(A_\theta-H_\theta I)|_{\Delta^\perp}$ are
$$
\delta_1=\phi_\theta\sqrt{(n-1)/2n}\;\;\mbox{and}\;\;\delta_2=-\phi_\theta\sqrt{(n-1)/2n}
$$
and, similarly, the eigenvalues of $(A-\lambda I)|_{\Delta^\perp}$ are
$$
\lambda_1=\phi_1\sqrt{(n-1)/2n}\;\;\mbox{and}\;\;\lambda_2=-\phi_1\sqrt{(n-1)/2n},
$$
where $\phi_1$ is given by \eqref{phi} with respect to $f$. On the other hand, 
since 
$$\det((A-\lambda I)J_\theta)=\det((A-\lambda I),$$
for  $\det J_\theta=1$, and both $(A-\lambda I)$ and $(A-\lambda I)J_\theta$ are traceless (see \eqref{traces}), it follows that $(A-\lambda I)$ and $(A-\lambda I)J_\theta$  have the same eigenvalues. This and \eqref{finalkey} imply that
$$
\phi_1^2=\varphi_\theta^2\phi_\theta^2,
$$
hence the Moebius metrics of $f$ and $f_\theta$ coincide. \qed

\vspace*{5ex}

\noindent Universidade de S\~ao Paulo\\
Instituto de Ci\^encias Matem\'aticas e de Computa\c c\~ao.\\
Av. Trabalhador S\~ao Carlense 400\\
13566-590 -- S\~ao Carlos\\
BRAZIL\\
\texttt{mibieta@icmc.usp.br} and \texttt{tojeiro@icmc.usp.br}

\newpage

\end{document}